\input amstex
\magnification=\magstep1
\baselineskip=13pt
\documentstyle{amsppt}
\vsize=8.7truein
\NoRunningHeads
\def\per{\operatorname{per}}
\def\EE{\bold{E}}
\def\PP{\bold{P}}

\def\xx{\bold{x}}
\def\rk{\operatorname{rank}}

\title Low Rank Approximations of Symmetric Polynomials and Asymptotic Counting 
of Contingency Tables  \endtitle
\author Alexander Barvinok \endauthor
\address Department of Mathematics, University of Michigan, Ann Arbor,
MI 48109-1109, USA \endaddress
\email barvinok$\@$umich.edu  \endemail
\date March 2005 \enddate
\thanks This research was partially supported by NSF Grant DMS 0400617.
\endthanks
\abstract We represent the number of $m \times n$ non-negative integer matrices (contingency tables) 
with prescribed row sums and column sums as the expected value of the permanent of a non-negative random matrix with exponentially distributed entries. We bound the variance of the obtained estimator, from which it follows that if the row and column sums are bounded by a constant fixed in advance, we get a polynomial time approximation scheme for counting contingency tables.
We show that the complete symmetric polynomial of a fixed degree in $n$ variables can be $\epsilon$-approximated coefficient-wise by a sum of powers of $O(\log n)$ linear forms, from which it 
follows that if the row sums (but not necessarily column sums) are bounded by a constant, there 
is a deterministic approximation algorithm of $m^{O(\log n)}$ complexity to compute 
the logarithmic asymptotic of the number of tables. \endabstract
\keywords contingency tables, approximation algorithms, symmetric functions, rank 
\endkeywords 
\subjclass 05A16, 68R05, 68W20, 15A15 \endsubjclass
\endtopmatter
\document
 
 \head 1. Introduction and main results \endhead 
 
 \subhead (1.1) Contingency tables \endsubhead
 Contingency tables are non-negative integer matrices with prescribed row and column sums,
 called marginals. The problem of computing the number of contingency tables with 
 given marginals has attracted a lot of attention recently, see \cite{DG95}, \cite{D+97}, 
 \cite{Mo02}, \cite{CD03}. The counting problem is motivated by applications to statistics, 
 combinatorics, representation theory, and is interesting in its own right, cf. \cite{DG95}.
 
 Let us consider non-negative integer $m \times n$ matrices
 with the row sums $r_1, \ldots, r_m$ and the column sums $c_1, \ldots, c_n$ such that
 $r_1 + \ldots + r_m = c_1 + \ldots + c_n = N$. If the 
 number $m$ of rows and the number $n$ of columns are fixed in advance, the number of 
 such matrices can be computed in polynomial time (that is, in time polynomial in $\log N$)
 since the problem reduces to counting integer points in a polytope in fixed dimension, see
 \cite{Ba94}. In fact, one does not need to apply the counting algorithm in full generality since the 
 polytope in question, the transportation polytope of non-negative matrices with prescribed 
 row and column sums, is either totally unimodular, or a straightforward ``combinatorial'' 
 degeneration of a totally unimodular polytope, a fact used in  \cite{DS03}.
 
 If one of the dimensions (for example, the number of columns) is allowed to grow, the 
 exact counting becomes difficult. As is shown in \cite{D+97}, exact counting is $\#P$-hard 
 already for $2 \times n$ matrices. The hardness result uses that the total sum $N$ can 
 become exponentially large in $n$ (so that $\log N$ is polynomial in $n$). If the number 
 $m$ of rows is fixed, a dynamic programming based algorithm computes the number of tables
 in time polynomial in $N$, thus resulting in a pseudo-polynomial algorithm, cf. \cite{CD03}.
 
On the other hand, Dyer, Kannan, and Mount \cite{D+97} have shown that 
if all the marginal $r_i$, $c_j$ are not too small ($r_i=\Omega(n^2 m)$ and $c_j=\Omega(m^2 n)$), then the Monte Carlo based approach allows one to approximate the number of contingency tables within a prescribed relative error $\epsilon >0$ in time polynomial in $m, n$ and $\epsilon^{-1}$. In this case,
the number of tables is well approximated by the volume of the corresponding polytope.
Subsequently, Morris improved the bounds to $r_i=\Omega(n^{3/2} m \log m)$ and 
$c_j=\Omega(m^{3/2} n \log n)$.  Combining the dynamic programming approach with the volume 
approximation idea, Cryan and Dyer \cite{CD03} obtained a randomized polynomial time 
approximation algorithm in the situation when the number of rows is fixed. This was later 
generalized in \cite{C+04}.

Thus the most difficult case is that with $N$ ``moderately large'' with respect to $m$ and $n$.

If both row sums $r_i$ and column sums $c_j$ are small, 
A. B\'ek\'essy, P. B\'ek\'essy, and Koml\'os \cite{B+72} proved the asymptotic formula
 $${N! \over r_1! \cdots r_m! c_1! \cdots c_n!} \exp\biggl\{{2 \over N^2} \sum_{i,j} {r_i \choose 2} 
 {c_j \choose 2}\biggr\} \tag1.1.1$$
  for the number of tables assuming that $N \longrightarrow +\infty$ while the 
 marginals remain bounded by a constant, fixed in advance: $r_i, c_j \leq \rho$.
 In \cite{B+72}, the authors proved that the relative error of this approximation is $O\left(N^{-1/2} \log N\right)$ and conjectured 
that it is $O\left( N^{-1} \right)$. Essentially, formula (1.1.1) counts contingency tables 
with entries not exceeding 2.

 Good and Crook \cite{GC77} make a heuristic argument that the formula should be valid for contingency tables under more general conditions of $r_i c_j/N$ being small. 
  
 If $m=n$ and $r_i=c_j=2$, an explicit generating function for the number of tables is
known, see Corollary 5.5.11 of \cite{St99}, which leads to a pseudo-polynomial algorithm
to compute the number of such tables exactly. 

Suppose now that we count every table $(d_{ij})$ with weight 
$$\prod_{ij} {1 \over d_{ij}!} \tag1.1.2$$ (the Fisher-Yates or the multiple hypergeometric statistics).
In this case, the weighted number of tables with row sums $r_1, \ldots, r_m$ and
column sums $c_1, \ldots, c_n$ is exactly equal to 
$${N! \over r_1! \cdots r_m! c_1! \cdots c_n!}. \tag1.1.3$$  
\subhead (1.2) Symmetric polynomials \endsubhead For a positive integer $r$, the 
{\it complete symmetric polynomial} $h_r$ of degree $r$ in $n$ variables $x_1, \ldots, x_n$ is the 
 sum of all distinct monomials 
 $$\xx^a=x_1^{\alpha_1} \cdots x_n^{\alpha_n} \quad \text{where} \quad 
 \sum_{i=1}^n \alpha_i =r \quad \text{and} \quad \alpha_i \geq 0 \quad \text{for} \quad 
 i=1, \ldots, n.$$
A well-known  and easy to prove result states that the number of $m \times n$ contingency tables
 with row 
sums $r_1, \ldots, r_m$ and column sums $c_1, \ldots, c_m$ is equal to the 
coefficient of the monomial $x_1^{c_1} \cdots x_n^{c_n}$ in the product 
$$h_{r_1}(\xx) \cdots h_{r_m}(\xx), \tag1.2.1$$
see, for example, Proposition 7.5.1 of \cite{St99}.
Similarly, if $e_r$ is an {\it elementary symmetric polynomial} of degree $r$ in $x_1, \ldots, x_n$
(that is, the sum of all square-free monomials of degree $r$), then the coefficient 
of the monomial $x_1^{c_1} \cdots x_n^{c_n}$ in the product 
$$e_{r_1}(\xx) \cdots e_{r_m}(\xx)$$ is 
the number of 0-1 matrices with the row sums $r_1, \ldots, r_m$ and 
the column sums $c_1, \ldots, c_n$, see Proposition 7.4.1 of \cite{St99}.

Let us ``approximate'' every polynomial $h_r$ in the product (1.2.1) by 
the power $(x_1 + \ldots  +x_n)^r$. The monomial expansion of the power contains
all the same monomials $\xx^a$ of degree $r$, only the coefficient of the monomial 
$x_1^{\alpha_1} \ldots x_n^{\alpha_n}$ is equal not to 1 but to $r!/ \alpha_1! \cdots \alpha_n!$.
Consequently, the coefficient of $x_1^{c_1}  \cdots x_n^{c_n}$ in the 
product
$$\split (x_1 + \ldots + x_n)^{r_1} \cdots (x_1 + \ldots + x_n)^{r_m}=&(x_1 + \cdots + x_n)^N,
\\ \text{where} \quad &N=r_1 + \ldots + r_m \endsplit$$
is equal to $r_1! \cdots   r_m!$ times the number of contingency tables with the row sums \break
$r_1, \ldots, r_m$ and the column sums $c_1, \ldots, c_n$, given that the 
weight of the table $(d_{ij})$ is the hypergeometric weight (1.1.2). On the other hand, 
this coefficient is equal to $N!/c_1! \cdots c_n!$, from which we deduce (1.1.3). 

As follows from formula (1.1.1), the Fisher-Yates statistics provides a reasonably good 
approximation to the uniform measure on contingency tables if the row and column sums 
are small. However, if only the row sums $r_i$ are small but column sums $c_j$ are allowed to be 
large (for example, if $m \gg n$), the approximation (1.1.1) is no longer valid. 

 In this paper, we present an algorithm for asymptotic computation of the number of contingency 
 tables where the row sums $r_i$ are small (and column sums $c_j$ are allowed to be large).
 Namely, for any $\epsilon>0$ and a positive integer $\rho$, fixed in advance, 
 we present an algorithm, which, given positive integers $r_1, \ldots, r_m \leq \rho$ and 
 positive integers $c_1, \ldots, c_n$, approximates the number of contingency tables with row 
 sums $r_1, \ldots, r_m$ and column sums $c_1, \ldots, c_n$ within a factor of $(1-\epsilon)^N$,
 where $N=r_1 + \ldots + r_m=c_1 + \ldots + c_n$. The algorithm has a quasi-polynomial complexity
 of $m^{O(\log n)}$. We present the algorithm in Section 3.
  The algorithm is based on the observation 
 that $n$-variate complete symmetric polynomials $h_r$ for small (fixed) $r$ can be approximated 
 by polynomials of  $O(\log n)$ {\it rank}. Namely, we prove the following result.
 \proclaim{(1.3) Theorem} Let us fix a positive integer $r$ and an $\epsilon>0$. Then there 
 exists a constant $\kappa=\kappa(r, \epsilon) >0$ with the following properties. For any 
 integer $n \geq 2$, there exist $k \leq \kappa \ln n$ linear forms 
 $\ell_i: {\Bbb R}^n \longrightarrow {\Bbb R}$ such that for the polynomial  
 $$\tilde{h}_r=\sum_{i=1}^k \ell_i^r(\xx)=\sum \Sb \alpha_1, \ldots, \alpha_n \geq 0 \\ 
 \alpha_1 + \ldots + \alpha_n =r \endSb \tilde{h}_{r,a} \xx^a,$$
 we have
 $$(1-\epsilon)^r \leq \tilde{h}_{r,a} \leq (1+\epsilon)^r$$
 for all non-negative integer vectors $a=(\alpha_1, \ldots, \alpha_n)$ with 
 $\alpha_1 + \ldots + \alpha_n=r$.
 \endproclaim
Moreover, we present a polynomial time algorithm to construct forms $\ell_i$.
Similar result holds for elementary symmetric functions $e_r$, which leads to a counting 
algorithm for 0-1 matrices. 
\proclaim{(1.4) Theorem} Let us fix a positive integer $r$ and an $\epsilon>0$. Then there 
 exists a constant $\kappa=\kappa(r, \epsilon) >0$ with the following properties. For any 
 integer $n \geq 2$, there exist $k \leq \kappa \ln n$ linear forms 
 $\ell_{ij}: {\Bbb R}^n \longrightarrow {\Bbb R}$ such that for the polynomial  
 $$\tilde{e}_r=\sum_{i=1}^k \prod_{j=1}^r \ell_{ij}(\xx)=\sum \Sb I=\{i_1, \ldots, i_r\} \\
 1 \leq i_1 < i_2 < \ldots < i_r \leq n  \endSb \tilde{e}_{r,I} x_{i_1} \cdots x_{i_r},$$
 we have
 $$(1-\epsilon)^r \leq \tilde{e}_{r,I} \leq (1+\epsilon)^r$$
 for all $r$-subsets $I \subset \{1, \ldots, n\}$.
\endproclaim

Let us fix a positive integer $k$. Let us ``approximate'' every polynomial $h_r$ in the product 
 (1.2.1) by a homogeneous polynomial of degree $r$ that is a product of polynomials $h_s$ with 
 $s \leq k$. Then the coefficients of $x_1^{c_1} \cdots x_n^{c_n}$ in the product  (1.2.1) enumerates
 contingency tables with weights ``interpolating'' between the Fisher-Yates statistics for 
 $k=1$ and the uniform measure on tables for $k \geq \max\{r_1, \ldots, r_m\}$. 
 
 As a by-product of our approach we express the number of contingency tables as the 
 expectation of the permanent of a random matrix. The {\it permanent} of an $N \times N$ matrix 
 $A=(a_{ij})$ is expressed by the formula
 $$\per A=\sum_{\sigma \in S_N} \prod_{i=1}^N a_{i \sigma(i)},$$
 where $\sigma$ ranges over the symmetric group $S_N$ of all permutations of 
 the set \break $\{1, \ldots, N\}$. Recently, Jerrum, Sinclair, and Vigoda constructed a 
 randomized polynomial time approximation scheme to compute the permanent of a 
 given non-negative matrix \cite{J+04}. As a corollary, they obtained a randomized 
 polynomial time approximation scheme to count 0-1 matrices with prescribed row and column
 sums. 
 
 Recall that a random variable $\gamma$ is {\it standard exponential} if 
 $$\PP(\gamma \geq t)=\cases e^{-t}  &\text{for\ }  t \geq 0 \\ 1 &\text{for \ } t < 0. \endcases$$ 
 We obtain the following result.
 \proclaim{(1.5) Theorem} Given positive integers $r_1, \ldots, r_m$ and $c_1, \ldots, c_n$ 
 such that $r_1 + \ldots + r_m=c_1 + \ldots + c_n=N$, let us consider the 
 $N \times N$ random matrix $A$ constructed as follows. We represent the set of 
 rows of $A$ as a disjoint union of $m$ subsets $R_1, \ldots, R_m$, where $|R_i|=r_i$ 
 for $i=1, \ldots, m$ and the set of columns of $A$ as a disjoint union of $n$ subsets 
 $C_1, \ldots, C_n$, where $|C_j|=c_j$. Thus $A$ is split into $m n$ blocks $R_i \times C_j$.
  We sample $mn$ independent standard exponential 
 random variables $\gamma_{ij}$, $i=1, \ldots, m$, $j=1, \ldots, n$ and fill the entries of 
 the block $R_i \times C_j$ by the copies of $\gamma_{ij}$.
 Let $\alpha =\per A$, so $\alpha$ is a function of the random variables $\gamma_{ij}$.
 
 Then
 \roster
 \item The number of $m \times n$ contingency tables with
 the row sums $r_1, \ldots, r_m$ and column sums 
 $c_1, \ldots, c_n$ is equal to 
 $${\EE \alpha \over r_1! \cdots r_m! c_1! \cdots c_n!};$$
 \item
 We have
 $${\EE \alpha^2 \over \EE^2 \alpha} \leq 2^{2N};$$
\item Suppose that $r_i, c_j \leq \rho$ for some number $\rho$ and all $i,j$. Then there exists
a constant $\kappa=\kappa(\rho)$, such that
$${\EE \alpha^2 \over \EE^2 \alpha} \leq \kappa.$$
One can choose $\kappa(\rho)=\exp\Bigl\{\rho^2 (2 \rho)! \Bigr\}$.
 \endroster  
 \endproclaim  
 We prove Theorem 1.5 in Section 4. Let us fix a number $0<p<1$, for example $p=2/3$.
 As follows by the Chebyshev inequality, if the row and column sums
 are bounded in advance, the average of $O(\epsilon^{-2})$ permanents of randomly generated
 $N \times N$ matrices, with probability at least $p$ approximates the number of contingency tables 
 within a relative error $\epsilon$.  In view of \cite{J+04},
  we obtain a polynomial time approximation
 algorithm for counting contingency tables when the row and column sums are 
 bounded by a constant, fixed in advance. 
 
 \subhead (1.6) Counting with weights \endsubhead A natural generalization of the 
 counting problem is counting with multiplicative weights. Given an $m \times n$ matrix 
 $W=(w_{ij})$ of weights, let us define the weight of an $m \times n$ non-negative integer 
 matrix $D=(d_{ij})$ as 
 $$\prod_{ij} w_{ij}^{d_{ij}}.$$
 For example, if $w_{ij} \in \{0,1\}$ then the weight of $D$ is 1 if and only if 
 $d_{ij} >0 $ implies $w_{ij}=1$. In this case, weighted counting implies counting matrices 
 with allowed entries $(i,j)$ for which $w_{ij}=1$.  Our results for asymptotic counting of contingency
 tables with small row sums admit generalization to counting with weights, provided the rank of 
 the weight matrix $W$ is fixed.  Similarly, Theorem 1.5 admits a straightforward generalization 
 for counting with weights: the $\gamma_{ij}$ entry of matrix $A$ needs to be multiplied by 
 $w_{ij}$. Part (2) also remains valid, although Part (3) does not.  Finally, we note that the weighted modification of the Fisher-Yates statistics can 
 be easily expressed as a permanent.
 \proclaim{(1.7) Theorem} Given positive integers $r_1, \ldots, r_m$ and $c_1, \ldots, c_n$ 
 such that $r_1 + \ldots + r_m=c_1 + \ldots + c_n=N$, and a non-negative $m \times n$ matrix 
 $W=(w_{ij})$, let us consider the $N \times N$ matrix $A$ constructed as follows.
 We represent the set of rows of $A$ as a disjoint union of $m$ subsets $R_1, \ldots, R_m$, where $|R_i|=r_i$ for $i=1, \ldots, m$ and the set of columns of $A$ as a disjoint union of $n$ subsets 
 $C_1, \ldots, C_n$, where $|C_j|=c_j$. Thus $A$ is split into $mn$ blocks $R_i \times C_j$.
 Let us fill the entries of the block $C_i \times R_j$ by $w_{ij}$.
 Then the total weight of $m \times n$ contingency tables with the row sums $r_1, \ldots, r_m$ 
 and column sums $c_1, \ldots, c_n$, where the table $D=(d_{ij})$ is counted with the 
 weight
 $$\prod_{ij} {w_{ij}^{d_{ij}} \over d_{ij}!},$$
 is equal to 
 $${\per A \over r_1! \cdots r_m! c_1! \cdots c_n!}.$$
 \endproclaim
 We prove Theorem 1.7 in Section 4.
 
 \head 2. Preliminaries: a scalar product in the space of polynomials \endhead
 
 We will use a certain scalar product in the space $V_n$ of real $n$-variate 
 polynomials. There are many ways to define it. The most straightforward way is to 
 define the scalar product of two monomials
 $$\langle \xx^a, \xx^b \rangle=\cases \alpha_1 ! \cdots \alpha_n! &\text{if\ } 
 a=b=(\alpha_1, \ldots, \alpha_n) \\ 0 &\text{otherwise.} \endcases$$
 A more formal way is to write 
 $$\langle f, g \rangle =f(\partial) g(\xx) \Big|_{\xx=(0, \ldots, 0),}$$
 where $f(\partial)$ is the differential operator
 $$f(\partial)=f\left( {\partial \over \partial x_1}, \ldots, {\partial \over \partial x_n} \right).$$
 The most invariant way is to consider the complex space ${\Bbb C}^n$, the Gaussian measure 
 $\nu_n$ there with the density
 $${1 \over \pi^n} e^{-\|z\|^2} \quad \text{where} \quad 
 \|z\|^2=|\zeta_1|^2 + \ldots + |\zeta_n|^2 \quad \text{for} \quad z=(\zeta_1, \ldots, \zeta_n),$$
 and let  
$$\langle f, g \rangle =\int_{{\Bbb C}^n} f(z) \overline{g(z)} \ d\nu_n.$$ 
 From this representation or otherwise, cf. \cite{Ba96}, it follows that the scalar product is invariant under orthogonal 
 transformations of the coordinates: if $U$ is an orthogonal transformation of ${\Bbb R}^n$ and 
 $f_1$ and $g_1$ are defined by $f_1(x)=f(Ux)$ and $g_1(x)=g(Ux)$ respectively, then 
 $\langle f, g \rangle = \langle f_1, g_1 \rangle$. Equivalently, for a linear transformation 
 $A: {\Bbb R}^n \longrightarrow {\Bbb R}^n$, let us define the polynomial $Af$ by 
 $$Af(\xx)=f(A^{\ast} \xx) \quad \text{for} \quad \xx \in {\Bbb R}^n,$$
 where $A^{\ast}$ is the conjugate transformation. Then
 $$\langle Af, g \rangle =\langle f, A^{\ast} g \rangle.$$ 
 
 The importance of this scalar product for us is that we can express the coefficient of $\xx^a$ 
 in $f$ as the scalar product 
 $${\langle f, \xx^a \rangle \over \alpha_1 ! \cdots \alpha_n!}.$$
 
 \subhead (2.1) Complexity of computing the scalar product \endsubhead
 Suppose that $f$ and $g$ are $n$-variate homogeneous polynomials of degree $k$ given by their 
 monomial expansions
 $$f(\xx)=\sum_{a} f_{a} \xx^{a} \quad \text{and} \quad g(\xx)=\sum_{a} g_{a}
 \xx^{a}.$$
 Then, to compute $\langle f, g \rangle$ one needs to sum up at most ${n+k-1 \choose k}$ 
 terms:
 $$\langle f, g \rangle =\sum \Sb a=(\alpha_1, \ldots, \alpha_n) \endSb 
 \alpha_1 ! \cdots \alpha_n ! f_{\alpha} g_{a}$$
 Taking into account computation of factorials, one can compute the scalar product 
 using $O\left(k {n+k-1 \choose k} \right)$ arithmetic operations. In particular, if the number of variables $n$ is fixed, we get a polynomial time algorithm.
 We will also be interested in the case of $n=O(\log k)$, in which case we get an algorithm of 
 a quasipolynomial $k^{O(\log k)}$ complexity.
 
 Generally, if the polynomials $f$ and $g$ are defined by their ``black boxes'', which, for any
 given $\xx=(x_1, \ldots, x_n)$ compute  the values $f(\xx)$ and $g(\xx)$, we can obtain 
 the monomial expansions of $f(\xx)$ and $g(\xx)$ via the standard procedure of interpolation 
 in $O\left( {n+ k -1 \choose k}^3 \right)$ time (provided $n$ and $k$ are known in advance), 
 cf. \cite{KY91} for the sparse version.
 Again, if $n$ is fixed, we get a polynomial time algorithm and if $n=O(\log k)$, we get an 
 algorithm of a quasipolynomial complexity. 
 \bigskip
 The invariance of the scalar product under the action of the orthogonal group often allows us 
 to reduce the number of variables. 
 \subhead (2.2) The rank of a polynomial \endsubhead Let $f: {\Bbb R}^n \longrightarrow {\Bbb R}$ be a  
 polynomial. We say that $\rk f \leq r$ if there are $r$ linear 
 forms $\ell_i: {\Bbb R}^n \longrightarrow {\Bbb R}$, $i=1, \ldots, r$ and a polynomial
 $q: {\Bbb R}^r \longrightarrow {\Bbb R}$ such that 
 $$f(\xx)=q\left(\ell_1(\xx), \ldots, \ell_r(\xx) \right) \quad \text{for} \quad \xx=(x_1, \ldots, x_n).$$
 
 Suppose we want to compute the scalar product $\langle f, g \rangle$, where $\rk f \leq r$ 
 and $f$ is represented as a polynomial $q$ in linear forms $\ell_1, \ldots, \ell_r$.
 Let $e_1, \ldots, e_r$ be the 
 coordinate linear forms
 $$e_i(x_1, \ldots, x_n)=x_i \quad \text{for} \quad i=1, \ldots, r.$$
 Let $A$ be a linear transformation such that $A e_i =\ell_i$ for $i=1, \ldots, r$.
 Then 
 $$\langle f, g \rangle  =\langle q\left(Ae_1, \ldots, Ae_r\right), g \rangle =
 \langle Aq(e_1, \ldots, e_r), g \rangle =\langle q(e_1, \ldots, e_r), A^{\ast} g \rangle.$$
 Now we observe that $q(e_1, \ldots, e_r)$ is a polynomial in the first $r$ variables $x_1, \ldots, x_r$.
 Replacing $A^{\ast}g$ by the ``truncated'' polynomial $\hat{g}$ obtained from $A^{\ast} g$ 
 by setting $x_{r+1} = \ldots = x_n=0$, we reduce computation of $\langle f, g \rangle$ to computation
 of the scalar product of two $r$-variate polynomials
 $$\langle f, g \rangle = \langle q, \hat{g} \rangle.$$
 
 In practical terms, if the linear forms $\ell_i$ are defined by
 $$\ell_i(\xx)=\alpha_{i1} x_1 + \ldots + \alpha_{in} x_n,$$
 we fill the $n \times n$ matrix $A^{\ast}=(a_{ij})$ by letting $a_{ij}=\alpha_{ij}$ for $i \leq r$
  and arbitrarily for 
 larger $i$. Then we transpose $A^{\ast}$ to get $A$ and compute $\hat{g}(x_1, \ldots, x_r)$ by 
 substituting $x_{r+1} =\ldots = x_n=0$ into $g(A\xx)$, where $\xx$ is interpreted as the 
 $n$-column of variables $x_1, \ldots, x_n$. 
 
 We will also need the following result, which can be considered as a complex version of the 
 Wick formula, see for example,  \cite{Zv97}. Since the author was unable to locate it 
 in the literature, we present its proof here. 
 \proclaim{(2.3) Lemma} Let $f_i, g_i: {\Bbb R}^n \longrightarrow {\Bbb R}$, $i=1, \ldots, m$ 
 be linear forms and let $F=f_1 \cdots f_m$ and $G=g_1 \cdots g_m$ be their products.
 Let us define an $m \times m$ matrix $B=(b_{ij})$ by 
 $b_{ij} =\langle f_i, g_j \rangle$ for $i,j=1, \ldots, m$. 
 Then 
 $$\langle F, G \rangle =\per B.$$
 \endproclaim
 \demo{Proof} First, we establish the formula in the particular case when 
 $g_1 = \ldots = g_m =e_1$, the 1st coordinate linear form. 
 In this case $G=x_1^m$, so letting $u=(1, 0, \ldots, 0)$, we can write 
 $$\langle F, G \rangle =m! F(u)= m! f_1(u) \cdots f_m(u).$$
 On the other hand, $b_{ij}=f_i(u)$, so
 $$\per B=m! f_1(u) \cdots f_m(u).$$
 
 Next, we establish the formula when $g_1 = \ldots = g_m$. In this case, we can write 
 $g_i = A e_1$ for some linear transformation $A$ of ${\Bbb R}^n$.
 Hence 
 $$\langle F, G \rangle =\langle F, (A e_1)^m \rangle = \langle A^{\ast} F, e_1^m \rangle=
 \langle (A^{\ast} f_1) \cdots (A^{\ast} f_m), e_1^m \rangle.$$
 Then $A^{\ast} f_i$ are linear forms and as we already established, the scalar product 
 is equal to the permanent of the matrix with the entries 
 $$\langle A^{\ast} f_i, e_1 \rangle = \langle f_i, A e_1 \rangle =\langle f_i, g_j \rangle =b_{ij}.$$
 Finally, we establish the general case of the formula. Let us fix the forms $f_1, \ldots, f_m$ 
 and consider both $\langle F, G \rangle$ and $\per B$ as functions of the forms 
 $g_1, \ldots, g_m$. We observe that both $\langle F, G \rangle$ and $\per B$ are 
 multilinear and symmetric in $g_1, \ldots, g_m$. Hence we obtain the general case by 
 polarization. Namely, let us fix $g_1, \ldots, g_m$. For real variables
 $t=(\tau_1, \ldots, \tau_m)$, let us define the linear form $g_t=\tau_1 g_1 + \ldots + \tau_m g_m$.
 Let $G_t =g_t^m$ and let $B(t)$ be defined by $b_{ij}(t)=\langle f_i, g_t \rangle$.
 Then both $\langle F, G_t \rangle$ and  $\per B(t)$ are homogeneous polynomials of degree
 $m$ in $\tau_1, \ldots, \tau_m$. Moreover, since both $\langle F, G \rangle$ and $\per B$ 
 are multilinear and symmetric in $g_1, \ldots, g_m$, the coefficient  
 of $\tau_1 \cdots \tau_m$ in $\langle F, G_t \rangle$ is equal to $m! \langle F, G \rangle$ while the coefficient of 
 $\tau_1 \cdots \tau_m$  in $\per B(t)$ is equal to $m!\per B$. Since we already proved 
 that 
 $\langle F, G_t \rangle =\per B(t)$, the result follows.
 {\hfill \hfill \hfill \hfill} \qed 
 \enddemo
 
 \head 3. Low rank approximations of symmetric polynomials  \endhead
 
Let $\gamma$ be a random variable with the standard exponential distribution
$$\PP(\gamma > \tau)=\cases e^{-\tau} &\text{for\ } \tau\geq 0 \\  1 &\text{for\ } \tau< 0. \endcases$$
 Hence for all integer $\alpha \geq 0$,
 $$\EE \gamma^{\alpha} =\alpha!.$$
 We will use the following straightforward result.
 \proclaim{(3.1) Lemma} Let $\gamma_1, \ldots, \gamma_n$ be independent random variables
 having the standard exponential distribution.
 Then, for any $r \geq 0$,
 $$\EE (\gamma_1 x_1 + \ldots + \gamma_n x_n)^r=r! h_r(x_1, \ldots, x_n),$$
 the complete symmetric polynomial of degree $r$.
 \endproclaim
 \demo{Proof} We have 
 $$(\gamma_1 x_1 + \ldots + \gamma_n x_n)^r =\sum \Sb \alpha_1, \ldots, \alpha_n \geq 0 \\
 \alpha_1 + \ldots + \alpha_n=r \endSb {r! \over \alpha_1! \cdots \alpha_n!} 
 \gamma_1^{\alpha_1} \cdots \gamma_n^{\alpha_n} x_1^{\alpha_1} \ldots x_n^{\alpha_n}.$$
 Since $\EE \gamma_i^{\alpha_i} =\alpha_i!$, the proof follows.
 {\hfill \hfill \hfill} \qed
 \enddemo
 
In what follows, $\kappa=\kappa(r, \epsilon)$ may denote various constants depending on $r$ and 
$\epsilon$ only.

Given a ``treshold'' $\kappa>0$, we define the truncated random exponential variable by
$$\overline{\gamma}= \cases \gamma &\text{if\ } \gamma \leq \kappa \\ 0 &\text{if\ } \gamma > \kappa, \endcases$$
where $\gamma$ is the standard exponential random variable.
The following is straightforward.
\proclaim{(3.2) Lemma} Given $r$ and $\delta>0$, there exists a constant 
$\kappa=\kappa(r, \delta)$ such that for the truncated random variable $\overline{\gamma}$, one has
$$(1-\delta)  \alpha! \leq \EE {\overline{\gamma}}^{\alpha} \leq \alpha! \quad \text{for} \quad 
\alpha=0, \ldots, r.$$ 
\endproclaim
 {\hfill \hfill \hfill} \qed
 
 Simple estimates show that one can choose 
 $$\kappa=O\left(r \ln r + \ln {1 \over \delta}\right).$$
 
 Next, we are going to use a concentration inequality (Azuma's inequality) for
 the sum of independent bounded random variables, see, for example, Theorem A.16 of \cite{AS92}.

\proclaim{(3.3) Proposition} Let $\xi_1, \ldots, \xi_m$ be independent random variables
such that $\EE \xi_i =\beta$ for $i=1, \ldots, m$ and $|\xi_i - \beta| \leq \kappa$ for some constant $\kappa$.
Then, for all $\delta > 0$, 
$$\PP\Bigl\{ \Big| {\xi_1 + \ldots + \xi_m \over m} -\beta \Big|> \delta \Bigr\} \leq 2 e^{-m \delta^2/ 2\kappa^2}.$$
\endproclaim
 An important consequence of Proposition 3.3
 is that for $\delta$, $\kappa$, and $r$ fixed, we can make the bound 
 $2e^{-m \delta^2/2\kappa^2}$ less than $n^{-r}$ by choosing $m=O(\log n)$.
 
 Now we can prove Theorem 1.3. 
 \demo{Proof of Theorem 1.3} We choose a $\delta >0$ so that 
 $(1-\delta) \geq (1-\epsilon)^{1/2}$ and a threshold $\kappa=\kappa(r, \delta)$ so as 
 to satisfy the conditions of Lemma 3.2. Then we sample 
 the coefficients of the linear forms $\ell_i: {\Bbb R}^n \longrightarrow {\Bbb R}, i=1, \ldots, m$ independently at random from the truncated standard
 exponential distribution.
 Let 
 $$\tilde{h}_r={1 \over r! m} \sum_{i=1}^m \ell_i^r.$$
 Then, each coefficient $\tilde{h}_{r,a}$ of the monomial $\xx^{a}$ is the average of $m$ independent random samples
 of the random variable
 $$\xi_a={\overline{\gamma_1}^{\alpha_1} \cdots \overline{\gamma_n}^{\alpha_n} \over
 \alpha_1! \cdots \alpha_n!}. $$ 
 Since $r$ is fixed, all random variables $\xi_{a}$ remain uniformly bounded by some constant 
 depending on $r$ and $\epsilon$ only. Moreover, $(1-\epsilon)^{r/2} \leq \EE \xi_{a} \leq 1$. 
 Since for a fixed $r$, the number of 
 ${n+ r -1 \choose r}$ of monomials $\xx^{a}$ of multidegree $a$ is bounded by a 
 polynomial in $n$, by Proposition 3.3 we can choose $m=O(\log n)$ so that for each $a$,
 the probability that the average of $\xi_{a}$ does not lie within 
 the interval $[(1-\epsilon)^r, (1+\epsilon)^r]$  
 does not exceed $\left(3{n+r-1 \choose r}\right)^{-1}$. Then, with probability at least 
 $2/3$, the average $\tilde{h}_r$ satisfies the conditions of Theorem 1.3.
 {\hfill \hfill \hfill} \qed 
 \enddemo
 
 We sketch the proof of Theorem 1.4 below.
 \demo{Sketch of proof of Theorem 1.4}  With a 
 surjective map $\omega: \{1, \ldots, n\} \longrightarrow \{1, \ldots, r\}$
 we associate a homogeneous polynomial $p_{\omega}$ of degree $r$ in $n$ 
 variables $x_1, \ldots, x_n$, which is the product of $r$ linear 
 forms in $\xx=(x_1, \ldots, x_n)$: 
 $$p_{\omega}(\xx)=\prod_{i=1}^r \sum \Sb j \in \{1, \ldots, n\} \\ \omega(j)=i \endSb x_j.$$  
If $\omega$ is sampled from the uniform distribution on the space of all surjective maps 
 $\{1, \ldots, n \} \longrightarrow \{1, \ldots, r\}$ then the expectation $\EE p_{\omega}$ is a 
 positive multiple of the elementary symmetric polynomial $e_k(\xx)$. Now we approximate 
 $\EE p_{\omega}$ by a sample average of $O(\log n)$ polynomials $p_{\omega}$.
 To sample $\omega$, it suffices to sample $\omega(i)$ independently for $i=1, \ldots, n$ 
 and accept the resulting map if it is surjective. The map fails to be surjective with 
 probability at most $r(1-1/r)^n$, which is negligible if $r$ is fixed and $n$ grows.  
 {\hfill \hfill \hfill} \qed
 \enddemo
 
 \subhead (3.4) Derandomization \endsubhead
 Proofs of Theorem 1.3 and 1.4 allow us to construct polynomials $\tilde{h}_r$ and 
 $\tilde{e}_r$ by averaging $O(\log n)$ polynomials that are built from linear functions 
 with independent random coefficients. A closer look reveals that the coefficients do not have 
 to be independent, but only $r$-wise independent (that is, every $r$ coefficients should be 
 independent). If $r$ is fixed in advance, one can use constructions of small (polynomial size) sample 
 spaces to simulate such random variables, cf. Section 2 of Chapter 15 of \cite{AS92} and 
 \cite{E+98}. This leads to polynomial time deterministic algorithms for construction of polynomials
 $\tilde{h}_r$ and $\tilde{e}_r$ in Theorems 1.3 and 1.4.

\subhead (3.5) Asymptotic  counting of contingency tables \endsubhead
Now we can come up with an algorithm for asymptotic counting of tables.
Let us fix an $\epsilon >0$ and a positive integer $\rho$. Suppose that $r_1, \ldots, r_m \leq \rho$.
We construct polynomials $\tilde{h}_r$ as in Theorem 1.3. The coefficient 
of $x^{c_1} \cdots x_n^{c_n}$ in the product 
$$H(\xx)=\tilde{h}_{r_1} \cdots \tilde{h}_{r_m}$$ 
up to a factor of $(1-\epsilon)^N$
is equal to the number of contingency tables with the row sums $r_1, \ldots, r_m$ 
and the column sums $c_1, \ldots, c_n$.
 Theorem 1.3 implies that the rank of $H$ is $O(\log n)$. Hence, applying 
the algorithm of Section 2.2, we compute the required coefficient in $m^{O(\log n)}$ time.

This construction allows some extensions and ramifications. 

First, it extends to counting with weights
(cf. Section 1.6) provided the rank of the weight matrix $W=(w_{ij})$ is fixed in advance.
To this end, we approximate the polynomial 
$h_{r_i}(w_{i1}x_1, \ldots, w_{in} x_n)$ by the sample average of $O(\log n)$ powers of 
linear forms $\ell^{r_i}_i(\xx)$ for $\ell_i(\xx)=\sum_{j=1}^n w_{ij} \gamma_{ij} x_j$, where 
$\gamma_{ij}$ are independent exponential random variables. If $\rk W$ is fixed in advance, the forms 
used in the approximation of $h_{r_i}$ span a subspace of $O(\log n)$ dimension.

Second, we can compute approximately various other expressions of the type 
$\langle Q(\xx), H(\xx) \rangle$. For example, let $C_k \subset \{1, \ldots, N\}$ for $k=1, \ldots, n$
be subsets of integers and let
$$Q(\xx)=\prod_{k=1}^n \sum_{c \in I_k} {x_k^c \over c!}.$$
Then, up to a factor of $(1-\epsilon)^N$, the value of $\langle Q(\xx), H(\xx) \rangle$ is 
equal to the number of contingency tables with the row sums $r_1, \ldots, r_m$ and 
all possible column sums $c_1, \ldots, c_n$ such that $c_k  \in C_k$ for $k=1, \ldots, n$.

Finally, using Theorem 1.4 instead of Theorem 1.3 we obtain asymptotic enumeration
algorithms for 0-1 matrices.

\head 4. The estimator for the number of tables \endhead

In this Section, we prove Theorems 1.5 and 1.7.

\demo{Proof of Theorem 1.5} 
Let us define an $n$-variate polynomial 
$$H(\xx)=\prod_{i=1}^m h_{r_i} (\xx) \quad \text{for} 
\quad \xx=(x_1, \ldots, x_n),$$
where $h_r(\xx)$ is the complete symmetric polynomial of degree $r$. Then the 
number of contingency tables with the row sums $r_1, \ldots, r_m$ and column sums 
$c_1, \ldots, c_n$ is equal to the coefficient of $\xx^C=x_1^{c_1} \cdots x_n^{c_n}$ in $H(\xx)$.
Using the scalar product of Section 2, we can write this number as 
$${\langle \xx^C, H(\xx) \rangle \over c_1! \cdots c_n!}.$$
Using Lemma 3.1, we express $H(\xx)$ as the expectation of a product of linear forms.
Namely, we define random linear forms $\ell_i$ by 
$$\ell_i(\xx)=\sum_{j=1}^n \gamma_{ij} x_j,$$
where $\gamma_{ij}$ are independent exponential random variables.
Then, by Lemma 3.1,
$$H(\xx)={1 \over r_1! \cdots r_m!} \prod_{i=1}^m \EE \ell_i^{r_i}(\xx)=
{1 \over r_1 ! \cdots r_m!} \EE \prod_{i=1}^m \ell_i^{r_i}(\xx).$$  
Let us denote $L(\xx)=\prod_{i=1}^m \ell_i^{r_i} (\xx)$. Hence the number of 
contingency tables can be written as 
$${\EE \langle L(\xx), \ \xx^C \rangle \over c_1! \cdots c_n! r_1! \cdots r_m!} .$$
Since both $L(\xx)$ and $\xx^C$ are products of linear forms, by Lemma 2.3 their scalar 
product evaluates by the permanent of the matrix of pairwise scalar products of linear 
forms $\ell_i(\xx)$ and $e_j(\xx)=x_j$, which is the matrix $A$.
This proves Part (1) of the theorem.

Let $S_N$ be the symmetric group of all permutations of the set $\{1, \ldots, N\}$. Denoting the 
entries of $A$ by $a_{ij}$, we get
$$\alpha=\sum_{\sigma \in S_N} \prod_{i=1}^N a_{i \sigma(i)} \quad \text{and} \quad 
\alpha^2 =\sum_{\phi, \psi \in S_N} \prod_{i=1}^N a_{i \phi(i)} a_{i \psi(i)}.$$
Therefore,
$$\split &\left(\EE \alpha\right)^2 =\sum_{\phi, \psi \in S_N} \left(\EE \prod_{i=1}^N a_{i\phi(i)} \right)
\left(\EE \prod_{i=1}^N a_{i \psi(i)} \right) \quad \text{and} \\
&\EE \alpha^2 =\sum_{\phi, \psi \in S_N}\EE \left( \prod_{i=1}^N a_{i \phi(i)} a_{i \psi(i)}\right).
\endsplit$$  
Hence we represented $\EE^2 \alpha$ and $\EE \alpha^2$ as a sum of $(N!)^2$ terms parameterized
by pairs of permutations $(\phi, \psi)$. 

To prove Part (2), we show that every term in the expansion of $\EE^2 \alpha$ is at most 
$2^{2N}$ times the corresponding term in the expansion of $\EE \alpha^2$. Indeed, each term in 
the expansion of $\EE^2 \alpha$ is the product of the type
$$\EE \left( \prod_{ij} \gamma_{ij}^{u_{ij}} \right) \EE \left( \prod_{ij} \gamma_{ij}^{v_{ij}} \right)
=\left(\prod_{ij} u_{ij}! \right)  \left(\prod_{ij} v_{ij}! \right),$$
where $u_{ij}$ and $v_{ij}$ are non-negative integers such that 
$$\sum_{ij} u_{ij} = \sum_{ij} v_{ij}=N.$$
The corresponding term in the expansion of $\EE \alpha^2$ is 
$$\EE \left( \prod_{ij} \gamma_{ij}^{u_{ij} + v_{ij}} \right)=\prod_{ij} (u_{ij}+v_{ij})!.$$
Hence the ratio is
$$\prod_{ij} {(u_{ij}+v_{ij})! \over u_{ij}! v_{ij}!} \leq \prod_{ij} 2^{u_{ij}+v_{ij}} =2^{2N},$$
which proves Part (2).

To prove Part (3), we notice that $\EE \alpha \geq N!$, since the expectation of every term
is at least 1.  Let us 
consider a particular term
$$t_{\phi \psi}=\prod_{i=1}^N a_{i \phi(i)} a_{j \psi(j)}.$$ 
We have $\EE t_{\phi \psi} >1$ if and only if some of the entries $a_{ij}$ in the product belong to the same block $R_i \times C_j$.
On the other hand, the maximum number of entries $a_{ij}$ which belong to the same block 
does not exceed $2 \rho$. Therefore, if the number of blocks with more than one entry is $k$,
$$\EE t_{\phi \psi} \leq \left((2 \rho)!\right)^k.$$   

Let us bound the number of terms $t_{\phi \psi}$ with $k$ entries belonging to the same block.
We can choose a permutation $\phi \in S_N$ in $N!$ ways and a subset $I \subset \{1, \ldots, N\}$ 
of $k$ indices in ${N \choose k}$ ways. For each entry $a_{i \phi(i)}$ with $i \in I$ we identify 
the block where $a_{i \phi(i)}$ belongs. Hence we get $k$ or fewer blocks since some 
of them may coincide. Now, for each $i \in I$ there are at most $\rho$ choices of $j \in \{1, \ldots, N\}$
and at most $\rho$ choices of $\psi(j)$ such that $a_{i\phi(i)}$ and $a_{j \psi(j)}$ belong to 
the same block as $\phi(i)$. After that, there are $(N-k)!$ choices for $\psi(i)$ for $i \notin I$. Hence 
$$\EE \alpha^2 \leq \sum_{k=0}^N {(N!)^2 \rho^{2k}  \over k!} \left((2 \rho)!\right)^{2k} \leq 
\left(N!\right)^2 \exp\Bigl\{ \rho^2 (2 \rho)! \Bigr\},$$
from which the proof of Part (3) follows.
{\hfill \hfill \hfill} \qed  
\enddemo
 The bound in Part (3) is probably non-optimal.
 \demo{Proof of Theorem 1.7} Let us define linear forms $\ell_i$ by
 $$\ell_i(\xx)=\sum_{j=1}^n w_{ij} x_i \quad \text{for} \quad i=1, \ldots, m.$$
 Let 
 $$L(\xx)=\prod_{i=1}^m \ell_i(\xx).$$
 As follows from the discussion of Section 1.2, the number of weighted tables can be expressed
 as the coefficient of $\xx^C=x_1^{c_1} \cdots x_n^{c_n}$ in the product $L(\xx)$ divided 
 by $r_1! \cdots r_m!$. Using the scalar product of Section 2, we write the number of 
 weighted tables as 
 $${ \langle L(\xx), \ \xx^C \rangle \over c_1! \cdots c_n! r_1! \cdots r_m!}.$$
Since both $L(\xx)$ and $\xx^C$ are products of linear forms, by Lemma 2.3 their scalar 
product evaluates by the permanent of the matrix of pairwise scalar products of linear 
forms $\ell_i(\xx)$ and $e_j(\xx)=x_j$, which is the matrix $A$.
{\hfill \hfill \hfill} \qed
 \enddemo
 
\head References \endhead

\widestnumber\key{MMMM}

\ref\key{AS92}
\by N. Alon and J.H. Spencer
\book The Probabilistic Method. With an appendix by Paul Erd\"os
\bookinfo Wiley-Interscience Series in Discrete Mathematics and Optimization
 \publ John Wiley $\&$ Sons, Inc.
 \publaddr New York
 \yr 1992
 \endref

\ref\key{Ba94}
\by A. Barvinok 
\paper A polynomial time algorithm for counting integral points in polyhedra when the dimension is
fixed 
\jour Math. Oper. Res. 
\vol 19 
\yr 1994
\pages 769--779
\endref

\ref\key{Ba96}
\by A. Barvinok 
\paper Two algorithmic results for the traveling salesman problem
\jour Math. Oper. Res.
\vol 21 
\yr 1996
\pages  65--84
\endref

\ref\key{B+72}
\by A. B\'ek\'essy, P.  B\'ek\'essy, and J. Koml\'os 
\paper Asymptotic enumeration of regular matrices 
\jour Studia Sci. Math. Hungar. 
\vol 7 
\yr 1972
\pages  343--353
\endref

\ref \key{CD03}
\by M. Cryan and M.  Dyer
\paper A polynomial-time algorithm to approximately count contingency tables when the number of rows is constant. Special issue on STOC 2002 (Montreal, QC)
\jour J. Comput. System Sci.
\vol 67 
\yr 2003 
\pages 291--310
\endref

\ref\key{C+04}
\by M. Cryan, M. Dyer, and D. Randall
\paper Approximately counting integral flows and cell-bounded contingency tables
\paperinfo preprint
\yr 2004
\endref

\ref\key{DG95}
\by P. Diaconis and A. Gangolli 
\paper Rectangular arrays with fixed margins 
\inbook Discrete Probability and Algorithms (Minneapolis, MN, 1993)
\pages 15--41 
\paperinfo IMA Vol. Math. Appl.
\vol 72 
\publ Springer
\publaddr New York
\yr 1995
\endref 

\ref\key{DS03}
\by J. De Loera and B.  Sturmfels 
\paper Algebraic unimodular counting. Algebraic and geometric methods in discrete optimization 
\jour Math. Program., Ser. B
\vol  96 
\yr 2003
\pages 183--203
\endref

\ref\key{D+97}
\by M. Dyer, R. Kannan, and J.  Mount 
\paper Sampling contingency tables 
\jour Random Structures Algorithms
\vol 10 
\yr 1997
\pages 487--506
\endref

\ref\key{E+98}
\by G. Even, O. Goldreich, M. Luby, N. Nisan, and B. Veli\u ckovi\'c
\paper Efficient approximation of product distributions
\jour Random Structures $\&$ Algorithms
\vol 13
\pages 1--16
\yr 1998
\endref

\ref\key{GC77}
\by I.J. Good and J.F. Crook
\paper The enumeration of arrays and a generalization related to contingency tables
\jour  Discrete Math.
\vol 19 
\yr 1977
\pages 23--45
\endref

\ref\key{J+04}
\by M. Jerrum, A. Sinclair, and E. Vigoda
\paper A polynomial-time approximation algorithm for the permanent of a matrix with 
non-negative entries
\jour Journal of the ACM
\vol 51
\yr 2004
\pages 671--697
\endref

\ref\key{KY91}
\by E. Kaltofen and L. Yagati
\paper Improved sparse multivariate polynomial interpolation algorithms
\inbook  Symbolic and algebraic computation (Rome, 1988)
\pages 467--474
\paperinfo Lecture Notes in Comput. Sci.
\vol 358
\publ Springer
\publaddr Berlin
\yr 1989
\endref

\ref\key{Mo02}
\by B. Morris 
\paper Improved bounds for sampling contingency tables 
\jour Random Structures Algorithms 
\vol 21 
\yr 2002
\pages 135--146
\endref

\ref\key{St99}
\by R.P. Stanley 
\book Enumerative Combinatorics. Vol. 2
\bookinfo Cambridge Studies in Advanced Mathematics
\vol  62 
\publ Cambridge University Press
\publaddr Cambridge
\yr 1999
\endref

\ref\key{Zv97}
\by A. Zvonkin
\paper Matrix integrals and map enumeration: an accessible introduction
\paperinfo Combinatorics and physics (Marseilles, 1995) 
\jour Math. Comput. Modelling 
\vol 26 
\yr 1997
\pages  281--304
\endref

 \enddocument